\documentclass[
12pt,reqno]{amsart}
\usepackage{verbatim}
\usepackage{amssymb}

\usepackage{enumerate}
\usepackage[active]{srcltx}
\numberwithin{equation}{section}

\usepackage{t1enc}
\usepackage[utf8x]{inputenc}

\newtheorem{theorem}{Theorem}[section]
\newtheorem{lemma}[theorem]{Lemma}
\newtheorem{corollary}[theorem]{Corollary}

\newtheorem{problem}[theorem]{Problem}
\newtheorem{proposition}[theorem]{Proposition}

\theoremstyle{definition}
\newtheorem{definition}[theorem]{Definition}

\theoremstyle{remark}

\newcommand{\itemprefix}{}
\makeatletter
\newcommand{\myitem}{%
\item\protected@edef\@currentlabel{\itemprefix\theenumi}%
}
\makeatother

\newcommand{\setm}{\setminus}

\newcommand{\subs}{\subset}
\newcommand{\dom}{\operatorname{dom}}

\newcommand{\<}{\left\langle}
\renewcommand{\>}{\right\rangle}

\author[I. Juh\'asz]{Istv\'an Juh\'asz}
\address      { Alfr\'ed Rényi Institute of Mathematics%
}
\email{juhasz@renyi.hu}

\author[J. van Mill]{Jan van Mill}
\address{University of Amsterdam}
\email{j.vanMill@uva.nl}

\author[L. Soukup]{Lajos Soukup}
\thanks
  {
   }
\address
      { Alfr{\'e}d R{\'e}nyi Institute of Mathematics
}
\email{soukup@renyi.hu}

\author[Z. Szentmikl\'ossy]{Zolt\'an Szentmikl\'ossy}
\address{E\"otv\"os University of Budapest}
\email{szentmiklossyz@gmail.com}
\thanks{The first, third and fourth authors were supported
by NKFIH grant no. K129211.}

\title{The double density spectrum of a topological space}

\begin{document}

\begin{abstract}
It is an interesting, maybe surprising, fact that different dense subspaces of even "nice"
topological spaces can have different densities. So, our aim here is to
investigate the set of densities of {\em all} dense subspaces of a topological space $X$
that we call the {\em double density spectrum} of $X$ and denote by $dd(X)$.

We improve a
result from \cite{BJ} by showing that $dd(X)$ is always $\omega$-closed (i.e. countably closed) if $X$ is Hausdorff.

We manage to give complete characterizations of the double density spectra of Hausdorff and of regular spaces
as follows.

Let $S$ be a non-empty set of infinite cardinals. Then
\begin{enumerate}
  \item  $S = dd(X)$ holds for a Hausdorff space $X$ iff $S$ is $\omega$-closed and $\sup S \le 2^{2^{\min S}}$;

  \smallskip

  \item  $S = dd(X)$ holds for a regular space $X$ iff $S$ is $\omega$-closed and $\sup S \le {2^{\min S}}$.
\end{enumerate}

We also prove a number of consistency results concerning the double density spectra of compact spaces.
For instance:

\begin{enumerate}[(i)]
  \item
If  $\kappa = cf(\kappa)$ embeds in $\mathcal{P}(\omega)/fin$ and $S$ is any set of uncountable regular
cardinals $< \kappa$ with $|S| < \min S$, then there is a compactum $C$ such that $\{\omega, \kappa\} \cup S \subs dd(C)$,
moreover $\lambda \notin dd(C)$ whenever $|S| + \omega < cf(\lambda) < \kappa$ and $cf(\lambda) \notin S$.

  \smallskip

  \item It is consistent to have a separable compactum $C$ such that $dd(C)$ is not $\omega_1$-closed.
  \end{enumerate}

\end{abstract}

\maketitle

\section{Introduction}
The density $d(X)$ of a topological space $X$, i.e. the minimum cardinality of a dense subspace of $X$,
is one of the most important and thoroughly studied topological cardinal functions. One important
feature of it is that it is not monotone, that is we may have $d(Y) > d(X)$ for a subspace,
even for a dense subspace, $Y$ of $X$.
The aim of this paper is to investigate for a given space $X$ the set of densities of all its dense subspaces,
that we call the {\em double density spectrum} of $X$.

For any space $X$ we shall denote by $\mathcal{D}(X)$ the family of all dense subspaces of $X$.
(We shall also use the notation $\mathcal{D}(\tau)$ instead of $\mathcal{D}(X)$ where $\tau$ is the topology of $X$.)
Thus $$dd(X) = \{d(D) : D \in \mathcal{D}(X)\}$$ will denote the double density spectrum of $X$.
We shall also use $\mathcal{N}(X)$ (or $\mathcal{N}(\tau)$) to denote the family of all {\em nowhere dense} subsets of $X$.

Of course, we have $d(X) = \min dd(X)$ and the double density spectrum had implicitly appeared  in the
definition of the cardinal function $\delta(X) = \sup dd(X)$ that was previously studied, in
chronological order, in \cite{WS}, \cite{BJ}, and \cite{JSh}. Clearly, we have $dd(X) \subs [d(X),\delta(X)]$
for any space $X$.

We are going to say that the space $X$ is $d$-{\em stable} if for every $U \in \tau^+(X)$
we have $d(U) = d(X)$, where $\tau^+(X)$ denotes the family of all non-empty open sets in $X$.
It is trivial that $$\{U \in \tau^+(X) : U \text{ is } d\text{-stable}\}$$ forms a $\pi$-base for $X$.

It is obvious that every $D \in \mathcal{D}(X)$ includes $I(X)$, the set of all isolated points of $X$.
Consequently, if $I(X) \in \mathcal{D}(X)$ then $$dd(X) = \{d(X)\} = \{|I(X)|\}.$$ Moreover, if
$I(X) \notin \mathcal{D}(X)$ and $Y = X \setm \overline{I(X)}$ then every member of $\mathcal{D}(X)$
has a dense subset of the form $I(X) \cup D$ with $D \in \mathcal{D}(Y)$, hence
$$dd(X) = \{|I(X)| + \kappa : \kappa \in dd(Y)\}.$$

Consequently, if we know the double density spectrum of all crowded spaces, we may easily compute the
the double density spectrum of any space. Therefore, in what follows, by {\em space} we shall always
mean a {\em crowded} Hausdorff space.

We end this section by presenting two propositions in which some very simple and basic but useful
properties of the double density spectra of spaces are collected.

\begin{proposition}\label{pr:1}
The following are satisfied for any space $X$.
\begin{enumerate}[{(a)}]
  \item If $Y \in \mathcal{D}(X)$ then $dd(Y) \subs dd(X)$; if even
  $Int(Y) \in \mathcal{D}(X)$, or equivalently: $X \setm Y  \in \mathcal{N}(X)$,
  then $dd(Y) = dd(X)$.

  \smallskip

  \item For every $U \in \tau^+(X)$ we have $dd(U) \setm d(X) \subs dd(X)$.

  \smallskip

  \item If $X$ is $d$-stable then $[X]^{< d(X)} \subs \mathcal{N}(X)$.
\end{enumerate}
\end{proposition}

\begin{proof}
The first part of
(a) is trivial and the second follows from $\mathcal{D}(Y) = \{D \cap Y : D \in \mathcal{D}(X)\}$.
(b) holds because for any $E \in \mathcal{D}(U)$ and $D \in \mathcal{D}(X)$ with
$|D| = d(X)$ we have $E \cup (D \setm \overline{U}) \in \mathcal{D}(X)$.
Finally, (c) is just a reformulation of what $d$-stable means.
\end{proof}

\section{Characterizing the double density spectra}

What kind of sets of (infinite) cardinals may occur as the double density spectrum of a space?
The aim of this section is to answer this natural question.

Since $|X| \le 2^{2^{d(X)}}$ holds for any Hausdorff space $X$ and $d(X) = \min dd(X)$,
it is obvious that $S = dd(X)$ implies $\sup S \le 2^{2^{\min S}}.$ Moreover, if $X$ is regular then
by $\delta(X) \le \pi(X) \le w(X) \le 2^{d(X)}$ we even have $\sup S \le 2^{\min S}.$

A set of cardinals is $\omega$-closed if it contains the supremum of all its countable subsets.
Our next theorem yields a less obvious necessary condition for the validity of $S = dd(X)$.

\begin{theorem}\label{tm:wcl}
The double density spectrum $dd(X)$ of any space $X$ is $\omega$-closed.
\end{theorem}

\begin{proof}
Assume that $$S = \{\kappa_n : n < \omega\} \in [dd(X)]^\omega,$$ then we have to show that $\kappa = \sup S \in dd(X)$.
We may, of course, assume that $\kappa \notin S$.

By part (b) of Proposition \ref{pr:1}, $\kappa \in dd(X)$ follows from $\kappa \in dd(U)$ for some $U \in \tau(X)$.
So, assume $\kappa \notin dd(U)$ for all $U \in \tau(X)$ and define for each $U \in \tau^+(X)$ $$\lambda(U) = \sup (\kappa \cap dd(U)).$$
Clearly, the family $$\mathcal{L} = \{U \in \tau^+(X) : \,\forall\, V \in \tau^+(U)\,(\lambda(V) = \lambda(U))\}$$
forms a $\pi$-base for $X$.

If there is $U \in \mathcal{L}$ such that $\lambda(U) = \kappa$ then we may take a disjoint collection $\{V_n : n < \omega\} \subs \tau^+(U)$
because $X$ is both crowded and Hausdorff. But for any $n < \omega$ then $\,\lambda(V_n) = \lambda(U) = \kappa$ implies
that there is $D_n \in \mathcal{D}(V_n)$ with $\kappa_n < |D_n| < \kappa$. Obviously, then $D = \bigcup \{D_n : n < \omega\}$
is a dense subset of $V = \bigcup \{V_n : n < \omega\}$ such that $d(D) = |D| = \kappa$, hence $\kappa \in dd(V)$.

The other possibility is that we have $\lambda(U) < \kappa$ for all $U \in \mathcal{L}$. In this case we take
a maximal disjoint subcollection $\mathcal{U}$ of $\mathcal{L}$ and note that $W = \cup \mathcal{U}$ is dense open in $X$.
Consequently, part (a) of Proposition \ref{pr:1} implies $S \subs dd(W) = dd(X)$.

We clearly also have $|\mathcal{U}| \le d(X) < \kappa$, hence $\sup \{\lambda(U) : U \in \mathcal{U}\} = \kappa$.
Indeed, if we had $\sup \{\lambda(U) : U \in \mathcal{U}\} = \mu < \kappa$ then we could choose $n < \omega$
such that $\mu \cdot |\mathcal{U}| < \kappa_n < \kappa$. But then there is $D \in \mathcal{D}(W)$ with $d(D) = \kappa_n$ and
for each $U \in \mathcal{U}$ we have $d(D \cap U) \le \mu$, which would imply $d(D) \le \mu \cdot |\mathcal{U}| < \kappa_n$, a contradiction.

Thus we may pick for all $n < \omega$ distinct, hence disjoint, $U_n \in \mathcal{U}$ so that $\kappa_n < \lambda(U_n) < \kappa$.
This again implies that there are $D_n \in \mathcal{D}(U_n)$ with $\kappa_n < |D_n| < \kappa$, hence
$D = \bigcup \{D_n : n < \omega\}$
is a dense subset of $U = \bigcup \{U_n : n < \omega\}$ such that $d(D) = |D| = \kappa$, similarly as above.
\end{proof}

It turns out that the necessary condition of $\omega$-closedness together with the
obvious cardinality restrictions mentioned above actually characterize the double density spectra
of Hausdorff and of regular spaces. Both characterization results will make use Cantor cubes.
We recall that the Cantor cube $\mathbb{C}_{\mu} = \{0,1\}^\mu$ of weight $\mu$ has density $\log \mu$, see e.g.
5.4 of \cite{J}. In particular, we have $d(\mathbb{C}_{2^\kappa}) = \log 2^\kappa \le \kappa$.

\begin{theorem}\label{tm:T2}
Let $S$ be an $\omega$-closed set of infinite cardinals such that $\min S = \kappa$ and $\sup S \le 2^{2^{\kappa}}$.
Then there is a  Hausdorff space $X$ with $S = dd(X)$.
\end{theorem}

\begin{proof}
The underlying set $X$ of our
promised space will be a dense subset of $\mathbb{C}_{2^\kappa}$ with a topology $\tau$ that is finer than the subspace
topology $\varrho$ on $X$ inherited from $\mathbb{C}_{2^\kappa}$. For later use we note that $\varrho$ is CCC if $X \in \mathcal{D}(\mathbb{C}_{2^\kappa})$.

To get $X$, we first fix pairwise disjoint subsets $\{X_\lambda : \lambda \in S\}$ of $\mathbb{C}_{2^\kappa}$ such that
$|X_\lambda| = \lambda$ and $X_\lambda$ is $\lambda$-dense in $\mathbb{C}_{2^\kappa}$ for each $\lambda \in S$.
(The latter means that $|U \cap X_\lambda| = \lambda$ for every non-empty open set in $\mathbb{C}_{2^\kappa}$.)
We can do this because $\mathbb{C}_{2^\kappa}$ may be partitioned into $2^{2^{\kappa}}$ dense subsets of size $\kappa$.
Indeed, the cosets of a $\kappa$-sized dense subgroup of $\mathbb{C}_{2^\kappa}$ form such a partition.
We then set $X  = \bigcup \{X_\lambda : \lambda \in S\}.$

On every $X_\lambda$ we consider the topology $\tau_\lambda$ generated by all sets of the form
$G \setm A$ where $G \in \varrho \upharpoonright X_\lambda$ and $A \in [X_\lambda]^{< \lambda}$.
We then let
$$\mathcal{N} = \{N \subs X : \forall\, \lambda \in S\,\big(N \cap X_\lambda \in \mathcal{N}(\tau_\lambda)\big)\}.$$
Clearly, $\mathcal{N}$ is a proper ideal on $X$.

Finally, our topology $\tau$ on $X$ is generated by all sets of the form $U \setm N$ where $U \in \varrho$
and $N \in \mathcal{N}$. Since $\emptyset \in \mathcal{N}$ we then have $\varrho \subs \tau$, hence
$\tau$ is Hausdorff. It is also obvious that $$\mathcal{N}(\tau \upharpoonright X_\lambda) = \mathcal{N}(\tau_\lambda) \supset [X_\lambda]^{< \lambda}$$
for all $\lambda \in S$.

It immediately follows from our definitions that $X_\lambda \in \mathcal{D}(\tau)$ and $d(X_\lambda,\tau) = \lambda$
for each $\lambda \in S$, hence we have $S \subs dd(X,\tau)$ and $$d(X,\tau) = \min\,dd(X, \tau) =  \min\,S.$$
It remains to show that $\mu \notin S$ implies $\mu \notin dd(X,\tau)$. As
$|X| = \sup\,S$, we may also assume that $d(X,\tau) < \mu \le \sup\,S$.

So, consider any $\tau$-dense set $D \in [X]^\mu$, we shall show that $d(D,\tau) < \mu$, hence $\mu \notin dd(X,\tau)$.
To see this, we first note that $\lambda \in S$ with $\lambda > \mu$ implies
$D \cap X_\lambda \in \mathcal{N}(\tau_\lambda)$, hence we have
$$D \cap \bigcup \{X_\lambda : \lambda \in S \setm \mu\} \in \mathcal{N}.$$
But then $D \cap \bigcup \{X_\lambda : \lambda \in S \setm \mu\}$ is $\tau$-nowhere dense,
hence the subset $$E = D \cap \bigcup \{X_\lambda : \lambda \in S \cap \mu\}$$ of $D$ is still $\tau$-dense.
So, if $|E| < \mu$ then we are done,
hence we may assume that $|E| = \mu$.

Next we note that for any
$U \in \varrho^+$ we have $E \cap U \notin \mathcal{N}$, hence  there is $\lambda \in S \cap \mu$ such that $E \cap U \cap X_\lambda$
is somewhere dense with respect to $\tau_\lambda$, i.e.
$$E \cap U \cap X_\lambda \notin \mathcal{N}(\tau_\lambda).$$
This in turn means that there is some $V \in \varrho^+$ with $V \subs U$ such that
$E \cap V \cap X_\lambda$ is $\tau_\lambda$-dense in $V \cap X_\lambda$.

Consequently, if $\mathcal{V}$ is a {\em maximal disjoint} collection of those $V \in \varrho^+$ for which
there is some $\lambda(V) \in S \cap \mu$ such that
$E \cap V \cap X_{\lambda(V)}$ is $\tau_{\lambda(V)}$-dense in $V \cap X_{\lambda(V)}$ then
$W = \cup\,\mathcal{V} \in \mathcal{D}(X,\varrho)$.

Since $\varrho$ is CCC, the family $\mathcal{V}$ is countable, hence
$$\lambda^* = \sup \{\lambda(V) : V \in \mathcal{V}\} \in S \cap \mu$$
because $S$ is $\omega$-closed and $\lambda^* \le \mu \notin S$.
Let us now put $$Y = \bigcup \{X_\lambda : \lambda \in S \text{ and } \lambda \le \lambda^* \}.$$
We claim that $E \cap Y$ is $\tau$-dense in $E$, and hence in $D$.

Indeed, assume that $U \setm N$ is any $\tau$-basic set where $U \in \varrho^+$ and $N \in \mathcal{N}$.
Then there is $V \in \mathcal{V}$ with $U \cap V \ne \emptyset$, hence
the choice of $\lambda(V)$ and $N \cap X_\lambda \in \mathcal{N(\tau_\lambda)}$
imply that $$E \cap U \cap V \cap X_{\lambda(V)} \setm N \ne \emptyset.$$
But $X_{\lambda(V)} \subs Y$ for all $V \in \mathcal{V}$, hence
$E \cap Y \cap (U \setm N) \ne \emptyset$ as well. Since $$|E \cap Y| \le |Y| = \lambda^* < \mu,$$
we thus have $d(D,\tau) \le  d(E,\tau) \le \lambda^* < \mu$, and
the proof is completed.
\end{proof}

\bigskip

The corresponding characterization of the double density spectra
of regular spaces will be immediate from our following result.

\begin{theorem}\label{tm:T3}
If $S$ is any $\omega$-closed set of infinite cardinals such that $\min S = \kappa$ and $\sup S \le 2^{\kappa}$
then there is a dense subspace $X$ of the Cantor cube $\mathbb{C}_{2^\kappa}$ of weight $2^\kappa$ with $S = dd(X)$.
\end{theorem}

\begin{proof}
We first claim that it suffices
to prove our theorem in the case in which $\log 2^\kappa = \kappa$.
Indeed, assume this and consider the case with $\log 2^\kappa < \kappa$. We may then apply the previous case
for the set $S' = \{\log 2^\kappa\} \cup S$ and obtain $X' \in \mathcal{D}(\mathbb{C}_{2^\kappa})$ such that
$S' = dd(X')$. Of course, we may also obtain $Y \in \mathcal{D}(\mathbb{C}_{2^\kappa})$ such that $dd(Y) = \{\kappa\}$,
for instance take $Y$ homeomorphic with  $\Sigma \times D$ where $\Sigma$ is the $\sigma$-product in $\mathbb{C}_{\kappa}$
and $D \in \mathcal{D}(\mathbb{C}_{2^\kappa})$ with $|D| \le \kappa$.
But then it is easy to check that $dd(X' \oplus Y) = S$, while it is also obvious that the topological sum $X' \oplus Y$
is homeomorphic to a dense subspace of $\mathbb{C}_{2^\kappa}$.

So, assume $d(\mathbb{C}_{2^\kappa}) = \log 2^\kappa = \kappa$ and note that this implies
$$[\mathbb{C}_{2^\kappa}]^{< \kappa} \subs \mathcal{N}(\mathbb{C}_{2^\kappa}).$$
We also fix $D \in \mathcal{D}(\mathbb{C}_{2^\kappa})$ with $|D| = \kappa$.

To prepare for the construction of our promised space $X$, we introduce the following definition.
For every set of indices $I \subs 2^\kappa$ we denote by $\sigma(I)$ the set of all functions
$f \in \{0,1\}^I$ with $|supp(f)| < \omega$, where $supp(f) = \{i \in I : f(i) = 1\}$.

Then we fix a disjoint family of sets $\{I_\lambda : \lambda \in S\}$ such that $I_\lambda \in [2^\kappa]^\lambda$,
and for each $\lambda \in S$ we let $K_\lambda = 2^\kappa \setm I_\lambda$. We then put
$$X_\lambda = \{x \in \mathbb{C}_{2^\kappa} : x \upharpoonright I_\lambda \in \sigma(I_\lambda) \text{ and }\exists d \in D (x \upharpoonright K_\lambda = d \upharpoonright K_\lambda)\}.$$
For every $\lambda \in S$ and $x \in X_\lambda$ we may then fix $d(x,\lambda) \in D$ such that
$x \upharpoonright K_\lambda = d(x,\lambda) \upharpoonright K_\lambda$.

Finally, we define the dense subspace $X$ of $\mathbb{C}_{2^\kappa}$ that we are looking for by $X = \bigcup \{X_\lambda : \lambda \in S\}$.
Clearly, we have $X_\lambda \in \mathcal{D}(\mathbb{C}_{2^\kappa})$ with $d(X_\lambda) = |X_\lambda| = \lambda$, moreover
$[X_\lambda]^{< \lambda} \subs \mathcal{N}(\mathbb{C}_{2^\kappa})$ for all $\lambda \in S$. This easily implies
$S \subs dd(X \cap U)$ for any $U \in \tau^+(\mathbb{C}_{2^\kappa})$. We claim that actually $S = dd(X \cap U)$ holds whenever $U \in \tau^+(\mathbb{C}_{2^\kappa})$.

The proof of this claim is indirect, so assume that $$\mu =  \min \bigcup \{dd(X \cap U) \setm S : U \in \tau^+(\mathbb{C}_{2^\kappa})\}$$
is well defined. Note that we then have $\mu > \kappa$. In what follows, we fix $U \in \tau^+(\mathbb{C}_{2^\kappa})$
with $\mu \in dd(X \cap U)$.
By definition, this means that there is $Y \in \mathcal{D}(X \cap U) \subs \mathcal{D}(U)$ with $d(Y) = |Y| = \mu$.

Let us recall that the $d$-stable open sets form a $\pi$-base in any space. We apply this to the space $Y$
and take a maximal disjoint collection $\mathcal{V} \subs \tau^+(U)$ such that
$Y \cap V$ is $d$-stable for each $V \in \mathcal{V}$.

We claim that there is a $V \in \mathcal{V}$ for which $d(Y \cap V) = \mu$. Indeed, if we have $d(Y \cap V) < \mu$ for some
$V \in \mathcal{V}$ then $d(Y \cap V) \in dd(X \cap V)$ and the minimality of $\mu$ together imply
$d(Y \cap V) \in S$. Consequently, as $\mathcal{V}$ is countable and
$S$ is $\omega$-closed, we cannot have $d(Y \cap V) < \mu$ for all $V \in \mathcal{V}$ because that would imply
$$d(Y) \le \sum \{d(Y \cap V) : V \in \mathcal{V}\} = \lambda \in S$$ with $\lambda < \mu$,
contradicting that $d(Y) = \mu$.

The family $\mathcal{E}$ of the elementary open sets $$[\varepsilon] = \{x \in \mathbb{C}_{2^\kappa} : \varepsilon \subs x\}$$
forms a base for $\mathbb{C}_{2^\kappa}$,
with $\varepsilon$ running through $Fn(2^\kappa, 2)$, the set of all finite partial functions from $2^\kappa$ to $2$. So, if we
fix $\varepsilon \in Fn(2^\kappa, 2)$ so that $[\varepsilon] \subs V \in \mathcal{V}$ with $d(Y \cap V) = \mu$
then $E = Y \cap [\varepsilon] \in \mathcal{D}([\varepsilon])$ is $d$-stable with $d(E) = |E| = \mu$.
This clearly implies $[E]^{< \mu} \subs \mathcal{N}(\mathbb{C}_{2^\kappa})$.

Consequently, for all $\lambda \in S$ we have $E \cap X_\lambda \in \mathcal{N}(\mathbb{C}_{2^\kappa})$.
Indeed, if $\lambda < \mu$ this follows from $|X_\lambda| = \lambda < \mu$, as we have just seen,
and if $\lambda > \mu$ then this follows from $|E| = \mu < \lambda$ and $[X_\lambda]^{< \lambda} \subs \mathcal{N}(\mathbb{C}_{2^\kappa})$.

For every point $d \in D$ let us now define $$J_d = \{\lambda \in S : \exists\,e \in E \cap X_\lambda\,(e \upharpoonright K_\lambda = d \upharpoonright K_\lambda)\},$$
moreover let $D_0 = \{d \in D : |J_d| \ge \omega\}$ and $D_1 = D \setm D_0$.

If $d \in D_0$ then we may pick distinct  cardinals $\{\lambda_n : n < \omega\} \subs J_d$ and points $e_n \in E \cap X_{\lambda_n}$
that "witness" $\lambda_n \in J_d$, i.e. $e_n \upharpoonright K_{\lambda_n} = d \upharpoonright K_{\lambda_n}$.
But for any basic neighborhood $[d \upharpoonright a] \in \mathcal{E}$ of $d$, where $a \in [2^\kappa]^{< \omega}$,
there are only finitely many $n$ with $a \cap I_{\lambda_n} \ne \emptyset$, hence we have $e_n \in [d \upharpoonright a]$
for all but finitely many $n < \omega$. In other words, this means that the sequence $\{e_n : n < \omega\} \subs E$
converges to $d$.

But this means that there is a subset $E_0 \subs E$ with $$|E_0| \le |D_0| \times \omega \le \kappa < \mu$$ such that
$D_0 \subs \overline{E_0}$. Consequently, $D_0$ is nowhere dense in $\mathbb{C}_{2^\kappa}$ because $E_0$ is.
This, in turn, implies that there is $\varepsilon_1 \in Fn(2^\kappa)$ such that $\varepsilon_1 \supset \varepsilon$
and $D \cap [\varepsilon_1] \subs D_1$.

For every $d \in D_1$ we have $|J_d| < \omega$ and for each $\lambda \in J_d$ we may pick $e(d,\lambda) \in E \cap X_\lambda$
such that $e(d,\lambda) \upharpoonright K_\lambda = d \upharpoonright K_\lambda$. We shall now show that the set
$$F = \{e(d,\lambda) : d \in D_1 \text{ and } \lambda \in J_d \} \subs E$$
is dense in $E \cap [\varepsilon_1]$.

Indeed, for any $\eta \in Fn(2^\kappa, 2)$ with $\eta \supset \varepsilon_1$ the set
$$Z_\eta = \bigcup \{E \cap X_\lambda : \dom(\eta) \cap I_\lambda \ne \emptyset\} \in \mathcal{N}(\mathbb{C}_{2^\kappa}).$$
So, we can pick $x \in E \cap [\eta] \setm Z_\eta$ and then $x \in X_\lambda$ implies $I_\lambda \cap \dom(\eta) = \emptyset$,
hence $d(x,\lambda) \in D \cap [\eta] \subs D \cap [\varepsilon_1] \subs D_1$. But $x$ is a witness for $\lambda \in J_{d(x,\lambda)}$,
so $f = e(d(x,\lambda),\lambda) \in F \cap [\eta]$ because
$$f \upharpoonright \dom(\eta) = d(x,\lambda) \upharpoonright \dom(\eta) = x \upharpoonright \dom(\eta).$$

But this yields the required contradiction that completes our indirect proof,
because $|F| \le \kappa$, and hence $F$ is nowhere dense.
\end{proof}

From this we immediately obtain the following characterization result.

\begin{corollary}
The following statements are equivalent for a non-empty set $S$ of infinite cardinals:
\begin{enumerate}[{(i)}]
\item $S$ is $\omega$-closed and $\sup S \le 2^{\min S}$.

\smallskip

\item There is a 0-dimensional CCC space $X$ with $S = dd(X)$.

\smallskip

\item There is a regular space $X$ with $S = dd(X)$.
\end{enumerate}
\end{corollary}

\section{On the double density spectra of compact spaces}

The aim of this section is to present what we know about the double density spectra of compact (Hausdorff) spaces.
Unfortunately, unlike for the classes of Hausdorff or regular spaces, we do not have any full
characterization in this case. However, we do know that the criteria for
regular spaces are not sufficient for the class of compacta. Indeed, the main result of \cite{JSh}
says that $\pi(X) = \max\,dd(X)$ holds for any compactum $X$. In particular, this implies that
the double density spectrum of a compact space always "admits a top".

Of course, as $\delta(X) \le \pi(X) \le w(X) \le 2^{d(X)}$ holds for a compactum $X$, there are only
problems if $2^{d(X)} > d(X)^+$, i.e. if the GCH fails at $\kappa = d(X)$. This leads us to the
following question.

\begin{problem}\label{pr:cpt}
Assume that $2^\kappa > \kappa^+$. Is there a compact space $X$ with $d(X) = \kappa$ and $\pi(X) > \kappa^+$
such that $\kappa^+ \notin dd(X)$,
or at least such that $dd(X) \ne [\kappa,\,\pi(X)]$?
\end{problem}

We do not have a complete answer to these questions but we shall present below interesting consistency results concerning them.
These results suggest that, at least consistently, we have a considerable amount
of freedom about the double density spectra of compact spaces.

We shall actually concentrate on the perhaps most interesting case of $\kappa = \omega$, i.e.
the double density spectra of {\em separable} compact spaces. To do that,
we shall develop a general method of constructing $\sigma$-centered and 0-dimensional spaces on
$[\omega]^\omega$ from appropriate ideals on $\omega$. The required separable compacta
will be just the compactifications of these "ideal spaces".

So, in what follows, "ideal" will always mean an ideal on $\omega$ that contains all finite subsets of $\omega$.

\begin{definition}
For any finite sequence $s \in 2^{< \omega}$ and $I \subs \omega$ we let
$$B(s, I) = \{A \in [\omega]^\omega : s \subs \chi_A \text{ and } A \cap I \subs |s|\},$$
where $\chi_A$ is the characteristic function of $A$ and $|s|$ is the length of $s$.
\end{definition}

Now, fix  an ideal $\mathcal{I}$ on $\omega$, then we let $\tau_\mathcal{I}$ denote the topology
on $[\omega]^\omega$ generated by $$\mathcal{B}_\mathcal{I} = \{B(s,I) : s \in 2^{< \omega} \text{ and }  I \in \mathcal{I}\}.$$

It is not completely obvious but it is immediate from the next lemma that $\mathcal{B}_\mathcal{I}$
is actually a base for $\tau_\mathcal{I}$.

\begin{lemma}\label{lm:base}
If $B(s,I) \cap B(t,J) \ne \emptyset$ then  $$B(s,I) \cap B(t,J) = B(s \cup t, I \cup J).$$
\end{lemma}

\begin{proof}
Clearly, if $B(s,I) \cap B(t,J) \ne \emptyset$ then either $s \subs t$ or $t \subs s$. By symmetry,
we may assume that $s \subs t$, hence $s \cup t = t$. Now, it is again obvious that then
$$B(s,I) \cap B(t,J) = B(s,I) \cap B(t, I \cup J) = B(s,I) \cap B(t,I) \cap B(t,J).$$
From this we get that $B(s,I) \cap B(t,I) \ne \emptyset$ as well.
But for any $A \in B(s,I) \cap B(t,I)$ and $|s| \le k < |t|$ we then have $k \notin A \cap I$, hence $k \in I$
implies $\chi_A(k) = t(k) = 0$. But this implies
$B(t,I) \subs B(s,I)$, and so we can conclude that
$B(s,I) \cap B(t,J) = B(t, I \cup J).$
\end{proof}

Note that $\{B(s, \emptyset) : s \in 2^{< \omega}\}$ is just the standard base for the Baire space, hence $\tau_\mathcal{I}$
is Hausdorff, being finer than the Baire space topology.

It is also obvious that for every $s \in 2^{< \omega}$ the collection $\{B(s,I) : I \in \mathcal{I}\}$ is
centered because $\mathcal{I}$ is an ideal, hence $\tau_\mathcal{I}$ is $\sigma$-centered.

To see that $\tau_\mathcal{I}$ is 0-dimensional, we shall show that evey member of its base $\mathcal{B}_\mathcal{I}$
is also closed.
Indeed, if $A \notin B(s,I)$ then either $s \nsubseteq \chi_A$ and then for $t = \chi_A \upharpoonright |s|$ we have
$A \in B(t, \emptyset)$ and $B(t, \emptyset) \cap B(s,I) = \emptyset$, or $s \subset \chi_A$ and there is $k \in A \cap I \setm |s|$.
In the latter case let $t$ be any initial segment of $\chi_A$ with $|t| > k$, then again we have $B(t, \emptyset) \cap B(s,I) = \emptyset$.

It follows then that the space $X_\mathcal{I} = \<[\omega]^\omega,\tau_\mathcal{I}\>$ does have compactifications,
and every compactification $C$ of $X_\mathcal{I}$ is separable because $\tau_\mathcal{I}$ is $\sigma$-centered.
Thus we have $\min dd(C) = \omega$.

A subset $\mathcal{A} \subs [\omega]^\omega$ is called {\em full} if for every $A \in \mathcal{A}$ and $B =^* A$
we have $B \in \mathcal{A}$.  In other words, $\mathcal{A}$ is full iff it is the union of $=^*$-equivalence classes.

We claim that $\pi(X_\mathcal{I}) = cof(\mathcal{I})$. Indeed, this is because if
$\mathcal{A} \subs \mathcal{I}$ is full then clearly
$$\{B(s,I) : s \in 2^{< \omega} \text{ and } I \in \mathcal{A}\}$$ is a $\pi$-base for $X_\mathcal{I}$ iff
$\mathcal{A}$ is cofinal in $\mathcal{I}$.
If $C$ is any compactification of $X_\mathcal{I}$ we thus have  $$\pi(C) = \pi(X_\mathcal{I}) = cof(\mathcal{I}) = \max dd(C).$$

Next we are going to describe $dd(X_\mathcal{I})$;
this is of interest because for any compactification $C$ of $X_\mathcal{I}$ we have $dd(X_\mathcal{I}) \subs dd(C)$.
First we introduce some new terminology.

A set $\mathcal{A} \subs [\omega]^\omega$ is called {\em $\mathcal{I}$-avoiding} if for every  $I \in \mathcal{I}$
there is $A \in \mathcal{A}$ such that $A \cap I = \emptyset$. Clearly, if $\mathcal{A}$ is full then
$\mathcal{A}$ is $\mathcal{I}$-avoiding iff for every  $I \in \mathcal{I}$
there is $A \in \mathcal{A}$ such that $|A \cap I| < \omega$.

\begin{lemma}\label{lm:fdense}
If $\mathcal{A} \subs [\omega]^\omega$ is full then $\mathcal{A} \in \mathcal{D}(X_\mathcal{I})$ iff $\mathcal{A}$ is $\mathcal{I}$-avoiding.
Consequently, if there is an $\mathcal{I}$-avoiding full set $\mathcal{A} \subs [\omega]^\omega$
such that $|\mathcal{A}| = \lambda$ but no $\mathcal{B} \subs \mathcal{A}$ with $|\mathcal{B}| < \lambda$  is $\mathcal{I}$-avoiding
then $\lambda \in dd(X_\mathcal{I})$.
\end{lemma}

\begin{proof}
Clearly, if $\mathcal{A}$ is full then $\mathcal{A}$ is $\mathcal{I}$-avoiding iff $\mathcal{A} \cap B(s,I) \ne \emptyset$ for every
$B(s,I) \in \mathcal{B}_\mathcal{I}$, hence iff $\mathcal{A}$ is dense in $X_\mathcal{I}$.

The second part follows because for every infinite $\mathcal{B} \subs \mathcal{A}$ the full hull $\mathcal{B}^*$ of
$\mathcal{B}$ has the same cardinality as $\mathcal{B}$. Consequently, no $\mathcal{B} \subs \mathcal{A}$ with $|\mathcal{B}| < \lambda$  is
dense in $X_\mathcal{I}$ as $\mathcal{B}^*$ isn't.
\end{proof}

Our next task is to find conditions that will ensure $\lambda \notin dd(C)$ for some cardinal $\lambda$
and compactification $C$ of $X_\mathcal{I}$. We start with a definition.

\begin{definition}
A cardinal $\lambda$ is said to be a {\em strong caliber} of the space $X$
(in symbols: $\lambda \in scal(X)$) if for every $\mathcal{U} \in [\tau(X)]^\lambda$
there is  $\mathcal{V} \in [\mathcal{U}]^\lambda$ such that $Int\big(\bigcap \mathcal{V}\big) \ne \emptyset$.
\end{definition}

The following simple but very useful proposition yields a condition for $\lambda \notin dd(X)$.

\begin{proposition}\label{pr:scal}
If $cf(\lambda) \in scal(X)$ then $\lambda \notin dd(X)$.
\end{proposition}

\begin{proof}
We prove the contrapositive of our statement. So, assume that $\lambda \in dd(X)$ and
$D \in \mathcal{D}(X)$ is such that $|D| = \lambda$ and $[D]^{< \lambda} \cap \mathcal{D}(X) = \emptyset$.
Enumerate $D$ as $D = \{x_\alpha : \alpha < \lambda\}$
and let $F_\alpha = \overline{\{x_\beta : \beta < \alpha\}}$.
Then $U_\alpha = X \setm F_\alpha \in \tau^+(X)$ and it is clear that for any cofinal subset $S \subs \lambda$
we have $Int\big(\bigcap \{U_\alpha : \alpha \in S\}\big) = \emptyset$ because
$\bigcap \{U_\alpha : \alpha \in S\} \cap D = \emptyset$.
But this clearly implies  $cf(\lambda) \notin scal(X)$.
\end{proof}

Although the definition of $\lambda$  being a strong caliber of $X$ uses $\tau(X)$, the family of all open
sets, it can obviously be replaced with any $\pi$-base of $X$. This implies that if $X$ is
quasiregular, i.e. the family $RO(X)$ of all regular open sets in $X$ forms a $\pi$-base of $X$,
then for any $Y \in \mathcal{D}(X)$ we have $scal(Y) = scal(X)$. Indeed, this follows from the fact
that the map $U \mapsto U \cap Y$ is an isomorphism between $RO(X)$ and $RO(Y)$.
In particular, this means that for every compactification $C$ of $X_\mathcal{I}$ we have $scal(C) = scal(X_\mathcal{I})$.
Consequently, $cf(\lambda) \in scal(X_\mathcal{I})$ implies $\lambda \notin dd(C)$ by Proposition \ref{pr:scal}.

Now, it remains to find a condition on the ideal $\mathcal{I}$ that will imply $\lambda \in scal(X_\mathcal{I})$.
The following result yields us just such a condition. But first we need a new piece of terminology.
If $\lambda$ is a cardinal then we call the ideal $\mathcal{I}$ on $\omega$ {\em weakly $\lambda$-complete} if
for every $\mathcal{A} \in[\mathcal{I}]^{\lambda}$ there is $\mathcal{B} \in[\mathcal{A}]^{\lambda}$
with $\cup \mathcal{B} \in \mathcal{I}$.

\begin{lemma}\label{lm:idscal}
Let $\lambda$ be a cardinal with $cf(\lambda) > \omega$. If the ideal $\mathcal{I}$ on $\omega$ is
weakly $\lambda$-complete then $\lambda \in scal(X_\mathcal{I})$.
\end{lemma}

\begin{proof}
As noted above, it suffices to verify the requirements of $\lambda \in scal(X_\mathcal{I})$ for
members of the base $\mathcal{B}_\mathcal{I}$ of $X_\mathcal{I}$. Using $cf(\lambda) > \omega$,
however, any subfamily of $\mathcal{B}_\mathcal{I}$ of size $\lambda$ has a subfamily of the same
size of the form $\{B(s,I) : I \in \mathcal{A}\}$,
with all members having the same first coordinate $s$, moreover $\mathcal{A} \in[\mathcal{I}]^{\lambda}$.

By our assumption, we may then find $\mathcal{B} \in[\mathcal{A}]^{\lambda}$
with $\cup \mathcal{B} = J \in \mathcal{I}$.
But then for every
$A \in B(s,J)$ and  $I \in \mathcal{B}$ we have $A \cap I \subs A \cap J  \subs |s|$, hence $A \in B(s,I)$ as well.
Consequently we have $$B(s,J) \subs \bigcap \{B(s,I) : I \in \mathcal{B}\}.$$
\end{proof}

In the rest of this section we are going to construct ideals such that their ideal spaces and
their (separable) compactifications will provide us with a wide variety of
double density spectra. Of course, these will require that the continuum $\mathfrak{c}$
be very large.

Actually, we shall need a cardinal characteristic of the continuum to be large, namely
the smallest cardinal $\kappa$ that does not embed in $\mathcal{P}(\omega)/fin$,
which we shall denote by $\mathfrak{n}$. In other words, $\mathfrak{n}$ is
the smallest cardinal such that there is no mod finite strictly increasing sequence
of that length in $[\omega]^\omega$.

It is trivial that $\mathfrak{b}^+ \le \mathfrak{n} \le \mathfrak{c}^+$ and it is also
well-known that $MA$ implies $\mathfrak{n} = \mathfrak{c}^+$. On the other hand, Kunen proved in his
PhD thesis that if one adds any number of Cohen reals to a model of CH then $\mathfrak{n} = \mathfrak{b}^+ = \omega_2$
holds in the generic extension, see IV.7.53 of \cite{K}.
We shall obtain interesting ideals, and hence interesting separable compacta, when $\mathfrak{n}$ is large.

In what follows, it will be useful have the notation $\Re$ for the class of all uncountable regular cardinals.
Our next result is presented as a warm up, it is actually a very special case of a later result but, we think,
it is quite interesting in itself.

\begin{theorem}\label{tm:kappa}
If $\kappa \in \mathfrak{n} \cap \Re$ then there is a separable compactum $C$ of $\pi$-weight $\kappa$
such that $\kappa \cap \Re \subs scal(C)$. Consequently,
$\omega < \lambda < \kappa$ and $cf(\lambda) > \omega$ imply $\lambda \notin dd(C)$.
So, if $\kappa < \aleph_\omega$ then $dd(C) = \{\omega, \kappa\}$.
\end{theorem}

\begin{proof}
Assume that $\<A_\alpha : \alpha < \kappa\>$ is a mod finite strictly increasing sequence in $[\omega]^\omega$.
We may assume that $\bigcup \{A_\alpha : \alpha < \kappa\} = \omega$, hence if $\mathcal{I}$ is the ideal
generated by $\{A_\alpha : \alpha < \kappa\}$ then $[\omega]^{< \omega} \subs \mathcal{I}$.
We claim that if $C$ is any compactification of $X_\mathcal{I}$ then $C$ is as required.

Since $\kappa$ is regular, it is obvious that $\pi(C) = \pi(X_\mathcal{I}) = cof(\mathcal{I}) = \kappa$, hence
we have $\kappa \in dd(C)$. Actually, we even have $\kappa \in dd(X_\mathcal{I})$ because the full set
$$\mathcal{C} = \{\omega \setm A_\alpha : \alpha < \kappa\}^*$$
is clearly $\mathcal{I}$-avoiding such that no subset of $\mathcal{C}$ of smaller size is $\mathcal{I}$-avoiding.

Next we show that, for every $\lambda$ as above, $\mathcal{I}$ is weakly $\lambda$-complete. Indeed, $\lambda < \kappa = cf(\kappa)$
implies that for every $\mathcal{A} \in[\mathcal{I}]^{\lambda}$ there is some $\alpha < \kappa$ such that
$I \subs^* A_\alpha$ for all $I \in \mathcal{A}$. Then by $cf(\lambda) > \omega$ there is some $a \in [\omega]^{< \omega}$
such that for $\mathcal{B} = \{I \in \mathcal{A} : I \subs a \cup A_\alpha\}$ we have $|\mathcal{B}| = \lambda$.
Clearly, then $\cup \mathcal{B} \in \mathcal{I}$.
But then Lemma \ref{lm:idscal} implies $\lambda \in scal(X_\mathcal{I}) = scal(C)$, hence
$\lambda \notin dd(C)$ by Proposition \ref{pr:scal}.
\end{proof}

We do not know whether $\lambda \in dd(C)$ holds for $\lambda \in (\omega, \kappa)$ with $cf(\lambda) = \omega$.
In particular, what happens with $\aleph_\omega$ if, say, $\kappa = \aleph_{\omega+1}$?

However, and we mention it just as a curiosity, we do know that
$dd(X_\mathcal{I}) = \{\kappa\}$ if $\<A_\alpha : \alpha < \kappa\>$ is a tower.
Indeed, this means that for every $A \in [\omega]^\omega$ there is an $\alpha < \kappa$ with $|A \cap A_\alpha| = \omega$.
So, if $\mathcal{C} \subs [\omega]^\omega$ with $|\mathcal{C}| = |\mathcal{C}^*| < \kappa$ then
there is an $\alpha < \kappa$ with $|A \cap A_\alpha| = \omega$ for all $A \in \mathcal{C}^*$, hence by
Lemma \ref{lm:fdense} we have that $\mathcal{C}^*$
is not dense in $X_\mathcal{I}$ and so neither is $\mathcal{C}$.

\begin{theorem}\label{tm:S}
For every set $S \subs \mathfrak{n} \cap \Re$ there is a separable compactum $C$ such that

\begin{enumerate}[(i)]
\item $\pi(C) = \sup S$;

\smallskip

\item $S \subs dd(C)$;

\smallskip

\item if $\mu \in \Re \setm S$ with $\mu > |S|$ then $\mu \in scal(C)$, hence
$cf(\lambda) = \mu$ implies $\lambda \notin dd(C)$.
\end{enumerate}
\end{theorem}

We are going to present two quite different proofs of this result. Our first proof reduces it to Theorem
\ref{tm:kappa}.

\begin{proof}
We may apply Theorem \ref{tm:kappa} for every $\lambda \in S$ to obtain a separable compactum $C_\lambda$ of $\pi$-weight $\lambda$
such that $\lambda \cap \Re \subs scal(C_\lambda)$. We claim that the product $C = \prod \{C_\lambda : \lambda \in S\}$ is as required.
Since $|S| \le \mathfrak{c}$ and the product of at most $\mathfrak{c}$ separable spaces is separable, $C$
is a separable compactum. It is also clear that $\pi(C) = \sup S$.

We obtain (ii) because for every $\lambda_0 \in S$ we may apply Proposition \ref{pr:2} with $X = C_{\lambda_0}$ and
$Y = \prod \{C_\lambda : \lambda \in S \setm \{\lambda_0\}\}$ to conclude that $\lambda_0 \in dd(C)$.

Finally, to verify (iii), we consider $\mu \in \Re \setm S$ with $\mu > |S|$. If $\mu > \lambda \in S$ then $\pi(C_\lambda) = \lambda$
trivially implies  $\mu \in scal(C_\lambda)$,
while for $\mu <\lambda \in S$ we have $\mu \in scal(C_\lambda)$ by  the choice of $C_\lambda$.

But it is straightforward to check that if a regular cardinal $\mu$ is a strong caliber of every factor of
the product of fewer than $\mu$ spaces then it is also a strong caliber of the product as well.
Thus we indeed have $\mu \in scal(C)$.
\end{proof}

Our second proof produces a suitable ideal $\mathcal{I}$ on $\omega$ such that any compactification
$C$ of $X_\mathcal{I}$ is as required.

\begin{proof}
We shall actually produce an ideal $\mathcal{I}$ on $\omega$ with the following three properties:
\begin{enumerate}[(a)]
  \item $cof(\mathcal{I}) = \sup S$.

\smallskip

  \item For every $\lambda \in S$ there is an $\mathcal{I}$-avoiding full set $\mathcal{A} \subs [\omega]^\omega$
such that $|\mathcal{A}| = \lambda$ but no $\mathcal{B} \subs \mathcal{A}$ with $|\mathcal{B}| < \lambda$  is $\mathcal{I}$-avoiding.

\smallskip

  \item If $\mu \in \Re \setm S$ and $\mu > |S|$ then $\mathcal{I}$ is
weakly $\mu$-complete.
\end{enumerate}
It easily follows from our earlier results that then any compactification $C$ of $X_\mathcal{I}$ is as required.
We leave it to the reader to check the details of this, and we move to the definition of $\mathcal{I}$.

We first fix an almost disjoint family $\{Q_\lambda : \lambda \in S\} \subs [\omega]^\omega$; this is possible
because $|S| \le \sup S \le \mathfrak{c}$. Then, for each $\lambda \in S$, we pick a mod finite strictly increasing $\lambda$-sequence
$\{Q_{\lambda,\alpha} : \alpha < \lambda\} \subs [Q_\lambda]^\omega$. Then we define $\mathcal{I}$ as the ideal generated by
$\mathcal{Q} = \{Q_{\lambda,\alpha} : \lambda \in S,\, \alpha < \lambda\}.$
We may assume, without any loss of generality, that $\cup \mathcal{Q} = \omega$, hence $[\omega]^{< \omega} \subs \mathcal{I}$.

Then (a) holds trivially. Item (b) holds because for every $\lambda \in S$ the full set $\mathcal{Q_\lambda} = \{Q_\lambda \setm Q_{\lambda,\alpha} : \alpha < \lambda\}^*$
is $\mathcal{I}$-avoiding with
$|\mathcal{Q_\lambda}| = \lambda$ and no $\mathcal{R} \subs \mathcal{Q_\lambda}$ with $|\mathcal{R}| < \lambda$  is $\mathcal{I}$-avoiding.

To see (c), we note first that for every $I \in \mathcal{I}$ there is a finite set $a_I \subs S$ and a function $f_I$ with domain $a_I$ and
with $f_I(\lambda) \in \lambda$ for each $\lambda \in a_I$
such that $I \subs^* \cup \{Q_{\lambda, f_I(\lambda)} : \lambda \in a_I\}$.

Now, assume that $\mu \in \Re \setm S$ with $\mu > |S|$, moreover $\mathcal{A} \in [\mathcal{I}]^\mu$. By $\mu > |S|$ we may
then assume that for some fixed $a \in [S]^{< \omega}$ we have $a_I = a$ for all $I \in\mathcal{A}$. Next, we may choose
$\mathcal{B} \subs \mathcal{A}$ with $|\mathcal{B}| = \mu$ so that
for every $\lambda \in a \cap \mu$ there is a fixed $\alpha_\lambda < \lambda$ with $f_I(\lambda) = \alpha_\lambda$
whenever $I \in \mathcal{I}$.

But for every $\lambda \in a \setm \mu$, since $\lambda > \mu$ is regular, we can fix $\alpha_\lambda < \lambda$
such that $f_I(\lambda) < \alpha_\lambda$ for all $I \in \mathcal{B}$.
This clearly implies that $$I \subs^* J = \cup \{Q_{\lambda, \alpha_\lambda} : \lambda \in a\} \in \mathcal{I}$$ whenever $I \in\mathcal{B}$.
Now, then there is some $b \in [\omega]^{< \omega}$ such that $\mathcal{C} = \{I \in \mathcal{B} : I \subs J \cup b\}$ also
has size $\mu$, hence $\mathcal{I}$ is indeed weakly $\mu$-complete.
\end{proof}

The following immediate corollary of Theorems \ref{tm:S} and \ref{tm:wcl} is, we think, illuminating.

\begin{corollary}
Assume that $\mathfrak{n} > \aleph_\omega$. Then
\begin{enumerate}
  \item for every $a \in [\omega]^{< \omega}$ there is a separable compactum $C$ such that $dd(C) = \{\omega\} \cup \{\omega_n : n \in a\}$;

\smallskip

  \item for every $a \in [\omega]^{\omega}$ there is a separable compactum $C$ such that $dd(C) = \{\omega\} \cup \{\omega_n : n \in a\} \cup \{\aleph_\omega\}$.
\end{enumerate}
\end{corollary}

The above consistency results on the double density spectra of (separable) compact spaces need that
$\mathfrak{n}$ be large, i.e. that "long" well-ordered sequences embed into $\mathcal{P}(\omega)/fin$.
However, it is well-known that in some sense any partial order embeds into $\mathcal{P}(\omega)/fin$,
at least after passing to an appropriate generic extension of the ground model.
More precisely, for any poset $\mathbb{Q} = \<Q, \le\>$ there is a CCC notion of forcing forcing $\mathbb{P}$
such that $\mathbb{Q}$ embeds into $\mathcal{P}(\omega)/fin$ in the generic extension $V^\mathbb{P}$.
(We haven't found a direct reference to this folklore result, however it is an immediate consequence of items
2.5 and 2.7 of \cite{DP}.)

It is not a surprise then that more complicated posets embedded in $\mathcal{P}(\omega)/fin$
yield us further interesting consistency results on the double density spectra of (separable) compacta.
Our next result well illustrates this.

\begin{theorem}\label{tm:w1cl}
It is consistent to have a separable compact space $C$ such that $dd(C)$, which of course is $\omega$-closed,
is not $\omega_1$-closed.
\end{theorem}

\begin{proof}
We first let
$$S = \Re \cap \aleph_{\omega_1} \setm \{\omega_1\} = \{\omega_{\alpha + 1} : 1 \le \alpha < \omega_1\}.$$
Then we consider the poset $\mathbb{Q} = \<Q, \le\>$, where $Q = \Pi \{\lambda : \lambda \in S\}$ and for
$x, y \in Q$ we have $x \le y$ iff $x(\lambda) \le y(\lambda)$ for all $\lambda \in S$.
Next we move from our ground model $V$ to its generic extension $V^\mathbb{P}$, where $\mathbb{P}$
is a CCC notion of forcing forcing such that $\mathbb{Q}$ embeds into $\mathcal{P}(\omega)/fin$ in $V^\mathbb{P}$.
Warning: $Q = \big(\prod\{\lambda : \lambda \in S\}\big)^{V} \ne \big(\prod\{\lambda : \lambda \in S\}\big)^{V^\mathbb{P}}$!

From now on we work in $V^\mathbb{P}$, so we may fix a bijection $h$ of $Q$ into $[\omega]^\omega$ such that for any $x, y \in Q$
we have $x < y$ iff $h(x) \subs^* h(y)$. In what follows, we shall write $A_x$ instead of $h(x)$.
Not surprisingly, we let $\mathcal{I}$ be the ideal on $\omega$ generated by $\{A_x : x \in Q\}$.
We may assume, without any loss of generality, that $\cup \{A_x : x \in Q\} = \omega$, hence $[\omega]^{< \omega} \subs \mathcal{I}$.

We next show that $S \subs dd(X_\mathcal{I})$. Indeed, for every $\lambda \in S$ and $\alpha < \lambda$ we define
$x_{\lambda, \alpha} \in Q$ by putting $x_{\lambda, \alpha}(\lambda) = \alpha$ and $x_{\lambda, \alpha}(\mu) = 0$
for all $\mu \in S \setm \{\lambda\}$.

Let us then define $B_{\lambda, \alpha} = \omega \setm A_{x_{\lambda, \alpha} }$. Clearly, then
$\mathcal{B}_\lambda = \{ B_{\lambda, \alpha} : \alpha < \lambda\}^*$ is an $\mathcal{I}$-avoiding full set and we claim
that no smaller sized subset of it is $\mathcal{I}$-avoiding. This follows from the fact that, as $\mathbb{P}$ is CCC, every subset of $\lambda$
in $V^\mathbb{P}$ is covered by subset of $\lambda$ in $V$ of the same size, and hence is bounded in $\lambda$.
So, we have $\lambda \in dd(X_\mathcal{I})$ by Lemma \ref{lm:fdense}.

A very similar argument, using that $\mathbb{P}$ is CCC, yields that every subset $R$ of $Q$ in $V^\mathbb{P}$ with $|R| \le \omega_1$
is bounded in $\mathbb{Q}$, i.e. there is $y \in Q$ such that for every $x \in R$ we have $x \le y$. This, in turn, means that
we have $A_x \subs^* A_y$ for all $x \in R$, hence we can find a finite $a \subs \omega$ with $|\{x \in R : A_x \subs a \cup A_y\}| = \omega_1$ as well.
But this clearly implies that $\mathcal{I}$ is weakly $\omega_1$-complete, hence $\omega_1 \in scal(X_\mathcal{I})$.

Now, if $C$ is any compactification of $X_\mathcal{I}$ then, on one hand, we have $S \subs dd(C)$, and on the other $\omega_1 \in scal(C)$.
But the first item implies that $\aleph_{\omega_1}$ is an accumulation point of $dd(C)$, while $\aleph_{\omega_1} \notin dd(C)$ by the second.
Consequently, $dd(C)$ is indeed not $\omega_1$-closed.
\end{proof}

\end{document}